\documentclass[reqno,a4paper,twoside,11pt]{amsart}
\usepackage{subfigure}
\usepackage{amsthm,amsmath,amssymb}
\usepackage{eepic}
\usepackage{graphicx}
\usepackage{hyperref}
\usepackage{enumerate}
\theoremstyle{definition}
\newtheorem{definition}{Definition}[section]
\newtheorem{example}[definition]{Example}
\theoremstyle{remark}

\newtheorem{note}[definition]{Note}

\theoremstyle{plain}
\newtheorem{lemma}[definition]{Lemma}

\newtheorem{theorem}[definition]{Theorem}
\newtheorem{corollary}[definition]{Corollary}

\newcommand{\ddivides}[2]{{#1 \left\lvert\lvert {#2} \right.}}
\newcommand{\set}[1]{\left\{{#1}\right\}}
\def\setsuchas#1#2{\left\{\,{#1}\,\vrule\,{#2}\,\right\}}

\newcommand{\Nat}{\mathbb{N}}
\newcommand{\Primes}{\mathbb{P}}
\newcommand{\PPrimes}{\mathbb{PP}}
\newcommand{\C}{\mathbb{C}}
\newcommand{\R}{\mathbb{R}}

\begin{document}
\title[Simplicial complexes \& multiplicative functions]
{Simplicial complexes associated to certain subsets of natural
  numbers and its applications to multiplicative functions}
\author{Jan Snellman}
\address{Department of Mathematics\\
Stockholm University\\
SE-10691 Stockholm, Sweden}
\email{jans@matematik.su.se}
\maketitle

\begin{section}{Introduction}
  We call a set of positive integers closed under taking unitary
  divisors an \emph{unitary ideal}. It can
  be regarded as a simplicial complex. Moreover, a multiplicative
  arithmetical function
  on such a set corresponds to a function on the simplicial complex
  with the property that the value on a face is the product of the
  values at the vertices of that face. We use this observation to
  solve the following problems: 
  \begin{enumerate}[A.]
  \item Let \(r\) be a positive integer and \(c\) a real number.
    What is the maximum value that \(\sum_{s \in S}g(s)\) can obtain
    when \(S\) is a unitary ideal containing precisely \(r\) prime
    powers, and \(g\) is the multiplicative function determined by
    \(g(s)=c\) when \(s \in S\) is a prime power?
  \item Suppose that \(g\) is a multiplicative function which is \(\ge
    1\), and that we want to find the maximum of \(g(i)\) when \(1 \le i
    \le n\). At how many integers do we need to evaluate \(g\)?
  \item If \(S\) is a finite unitary ideal, and \(g\) is
  multiplicative and \(\ge 1\), then the maximum of \(g\) on \(S\)
  occurs at a facet, and any facet is optimal for some such \(g\). If
  \(W_1,\dots,W_\ell\) is an enumeration of the facets in some order,
  is there always a \(g\) as above so that
  \(g(W_1) \le g(W_2) \cdots \le g(W_\ell)\)?
  \end{enumerate}
\end{section}

\begin{section}{Unitary ideals and simplicial complexes}
  Let \(\Nat\) denote the non-negative integers and \(\Nat^+\) the
  positive integers, with subsets \(\Primes\) the prime numbers and
  \(\PPrimes\) the set of prime powers.
  Recall \cite{Afunc} that an \emph{unitary divisor} (or a \emph{block
  factor}) of 
   \(n \in \Nat^+\) is a divisor \(d\) such that
  \(\gcd(d,n/d)=1\). In this case, we write \(\ddivides{d}{n}\) or \(n=d \oplus
  n/d\). If \(\gcd(a,b)>1\) we put \(a \oplus b=0\).
  
  \begin{definition}
    A subset \(S \subset \Nat^+\) is a \emph{unitary ideal} if
    \begin{equation}
      \label{eq:1}
      s \in S, \quad d \in \Nat^+, \quad \ddivides{d}{s} 
      \qquad \implies \qquad d \in S
    \end{equation}
  \end{definition}
  
  \begin{definition}
    For any unitary ideal \(S \subset \Nat^+\) with \(X=X(S)=\PPrimes
    \cap S\), we define the simplicial complex \(\Delta(S)\) on
    \(X(S)\) by 
    \begin{equation}
      \label{eq:2}
      \begin{split}
      \Delta(S) \ni \sigma = \set{a_1,\dots,a_r} & \iff 
      a_1 a_2 \cdots a_r \in S \text{ and }\\
      & \quad \forall\, 1 \le i <j \le r:
      \gcd(a_i,a_j)=1 
      \end{split}
    \end{equation}
  \end{definition}
  Clearly, \(\Delta(S)\) is finite iff \(X(S)\) is finite iff \(S\) is
  finite. Furthermore:
  \begin{lemma}[\cite{Snellman:UniGenTrunc}]\label{lemma:anyfin}
    Any finite
  simplicial complex can be realized as \(\Delta(S)\) for some \(S\).
    \end{lemma}
    \begin{proof}
      Take as many  prime numbers as there are vertices in the
      simplicial complex, so that the vertex \(v_i\) corresponds to
      the prime number \(p_i\). For any
      \(\sigma=\set{v_{a_1},\dots,v_{a_r}}\) in the simplicial complex
      we let \(p_{a_1} p_{a_2} \cdots p_{a_r} \in S\).
    \end{proof}
    
    \begin{note}
      In what follows, we will sometimes regard elements in \(S\) as
      faces in \(\Delta(S)\), without explicitly pointing this
      out. We trust that the reader will not be confused by this.
    \end{note}
    
  Recall that an \emph{arithmetical function} is a function \(g:
  \Nat^+ \to \C\), and that an arithmetical function is
  \emph{multiplicative} iff 
  \begin{equation}
    \label{eq:3}
    g(ab) =g(a)g(b)
  \end{equation}
  whenever \(\gcd(a,b)=1\).
  Hence, a multiplicative function is determined by its values on
  \(\PPrimes\), and we have

  \begin{lemma}\label{lemma:simplexvalue}
    Let \(S\) be a unitary ideal, and \(g\) a multiplicative function.
    By abuse of notation, put \(g(\sigma) = g(a_1a_2
    \cdots a_r)\) if \(\sigma=\set{a_1,\dots,a_r} \in
    \Delta(S)\). Then \(g(\sigma)=g(a_1)g(a_2) \cdots g(a_r)\), so the
    value of \(g\) at a simplex is the product of the valueos of \(g\)
    at the  vertices of said simplex.
  \end{lemma}

  We'll be interested in three problems: 
  \begin{enumerate}
  \item Calculating the sum \(\sum_{s \in S} g(s)\),
  \item Maximizing \(g\) on \(S\),
  \item Finding the total orders on \(S\) induced by \(g\).
  \end{enumerate}
\end{section}

\begin{section}{Summing \(g\) on \(S\)}
  We henceforth assume that \(S\) is finite, with \(S \cap \PPrimes\)
  containing \(r\) elements, and that \(g\) is a
  multiplicative function. Put \(G(S)=\sum_{s \in S} g(s)\).

  Let us start with the simplest cases. If \(S\) consists entirely of
  prime powers then \(\Delta(S)\) consists of \(r\) isolated points,
  on which \(g\) can take any values. The other extreme is that \(s,t
  \in S\), \(\gcd(s,t)=1\) implies that \(st \in S\), and that all
  prime powers in \(S\) are in fact primes. Then \(S\)
  consists of all square-free products of these \(r\) primes, so
  \(\Delta(S)\) is an \((r-1)\)-dimensional simplex. In this case, it
  is easy to see that \(G(S)=\prod_{j=1}^r \bigl(1+g(p_j)\bigr)\), where
  \(p_1,\dots,p_r\) are the primes in \(S\).
  
  More generally, if \(\Delta(S)\) have \(\ell\) faces
  \(W_1,\dots,W_\ell\), with \(2^\ell < \lvert \Delta(S) \rvert\),
  then the following formula might be  useful. Put
  \begin{displaymath}
    \tilde{g}(\set{a_1,\dots,a_v}) = \prod_{j=1}^v (1+g(a_i))
    = \sum_{\sigma \in \set{a_1,\dots,a_v}} g(\sigma).
  \end{displaymath}
  The principle of Inclusion-Exclusion givesd
  \begin{equation}
    \label{eq:inclusionexclusion}
    G(S) = \sum_{i=1}^\ell \tilde{g}(W_i) \,-\, \sum_{1 \le i < j \le
        \ell} \tilde{g}(W_i \cap W_j) \,+\, \sum_{1 \le i < j < k \le
        \ell} \tilde{g}(W_i \cap W_j \cap W_k) \,-\, \dots
  \end{equation}

  If \(S\) is arbitrary, but \(g\) special in that it takes the same
  value on all prime powers, then \(G(S)\) is also easily calculable.
  \begin{lemma}\label{lemma:allsame}
    If there exists a \(c\) such that \(g(s)=c\) for all \(s \in S
    \cap \PPrimes\), then 
    \begin{equation}
      \label{eq:4}
      G(S) = (1,c,c^2,\dots,c^r) \cdot
    (1,f_0,f_1,\dots, f_{r-1}),
    \end{equation}
where \((f_0,f_1,\dots,f_{r-1})\) is
    the \emph{\(f\)-vector} of \(\Delta(S)\), i.e. \(f_i\) counts the
    number of \(i\)-dimensional (i.e. having \(i+1\) vertices)
    faces of \(\Delta(S)\). 
  \end{lemma}
  \begin{proof}
    A \((v-1)\)-dimensional simplex of \(\Delta(S)\) contributes
    \(c^v\) to \(G(S)\); there are \(f_{v-1}\) such simplexes, so the
    total contribution is \(c^v f_{v-1}\). Letting \(v\) range from
    \(0\) to \(r\) and summing yields the result.
  \end{proof}

  \begin{theorem}
    Let \(\Psi(r,c)\) denote the maximum that \(G(S)\) can obtain when
    \(\lvert S \cap \PPrimes \rvert = r\) and \(g(s)=c \in \R\) for
    all \(s \in S \cap \PPrimes\). Then, if \(r\) is odd,
    \begin{equation}
      \label{eq:5}
      \Psi(r,c) = 
      \begin{cases}
        1 + \sum_{i=1}^r c^i \binom{r}{i} & \text{ if } c > 0 \\
        1 + \sum_{i=1}^2 c^i \binom{r}{i} & \text{ if }
        \frac{-4}{n-3}< c < 0 \\
        1 + \sum_{i=1}^4 c^i \binom{r}{i} & \text{ if }
        \frac{-6}{n-5}< c < \frac{-4}{n-3}  \\
        1 + \sum_{i=1}^6 c^i \binom{r}{i} & \text{ if }
        \frac{-8}{n-7}< c < \frac{-6}{n-5}  \\
        \qquad \vdots & \qquad \vdots \\
        1 + \sum_{i=1}^{r-1} c^i \binom{r}{i} & \text{ if }
        c <\frac{-(n-1)}{2} 
      \end{cases}
    \end{equation}
    and if \(r\) is even
    \begin{equation}
      \label{eq:6}
      \Psi(r,c) = 
      \begin{cases}
        1 + \sum_{i=1}^r c^i \binom{r}{i} & \text{ if } c > 0 \\
        1 + \sum_{i=1}^2 c^i \binom{r}{i} & \text{ if }
        \frac{-4}{n-3}< c < 0 \\
        1 + \sum_{i=1}^4 c^i \binom{r}{i} & \text{ if }
        \frac{-6}{n-5}< c < \frac{-4}{n-3}  \\
        1 + \sum_{i=1}^6 c^i \binom{r}{i} & \text{ if }
        \frac{-8}{n-7}< c < \frac{-6}{n-5}  \\
        \qquad \vdots & \qquad \vdots \\
        1 + \sum_{i=1}^{r-2} c^i \binom{r}{i} & \text{ if }
        c <\frac{-(n-2)}{3} \\
        1 + \sum_{i=1}^{r} c^i \binom{r}{i} & \text{ if }
        c <-n
      \end{cases}
    \end{equation}
  \end{theorem}
  \begin{proof}
    Put \(\mathbf{c} =
    (c,c^2,\dots,c^r)\), \(\mathbf{f}=(f_0,f_1,\dots,f_{r-1}\).
    It follows from Lemma~\eqref{lemma:anyfin} that we must maximize
    \(\mathbf{f} \cdot \mathbf{c}\) over all possible \(f\)-vectors
    \(\mathbf{f}\) of simplicial complexes on \(r\) vertices. Since
    \(\mathbf{f} \cdot \mathbf{c}\)  is  a 
    linear function, it will suffice to evaluate \(\mathbf{f} \cdot
    \mathbf{c}\) 
    on a set
    of vertices that span the convex hull of \(f\)-vectors of
    simplicial complexes on \(r\) vertices. Kozlov \cite{Chull} showed
    that the set 
    \begin{equation}
      \label{eq:koz}
      \set{\tilde{F}_1,\dots,\tilde{F}_r}, \qquad \tilde{F}_i = 
      \Bigl( \binom{r}{1}, \binom{r}{2}, \dots, \binom{r}{i}, 0,
      \dots, 0 \Bigr)
    \end{equation}
    is minimal with the property that its convex hull contains all
    \(f\)-vectors of 
    simplicial complexes on \(r\) vertices. Hence, it is enough to
    decide which of the \(r\) numbers
    \begin{equation}
      \label{eq:kn}
      \begin{split}
        K_1  = \tilde{F}_1 \cdot \mathbf{c} &= cn\\
        K_2  = \tilde{F}_2 \cdot \mathbf{c} &= cn + c^2 \binom{n}{2}\\
        \vdots &  \\
        K_r  = \tilde{F}_r \cdot \mathbf{c} &= \sum_{i=1}^r c^i \binom{n}{i}
      \end{split}
    \end{equation}
    is the greatest.

    Clearly, if \(c>0\), then \(K_r\) is the greatest. If \(c<0\) and
    \(r\) is odd then
    we always have the  inequalities shown in Figure~\ref{fig:rodd}.
    \setlength{\unitlength}{0.6cm}
    \begin{figure}[hbtp]
      \tiny{
    \begin{picture}(20,3)
      \multiput(0,0)(4,0){4}{\circle{1}}
      \multiput(2,2)(4,0){4}{\circle{1}}
      \multiput(0.5,0.5)(4,0){4}{\line(1,1){1}}
      \multiput(3.5,0.5)(4,0){4}{\line(-1,1){1}}

      \put(18,2){\circle{1}}
      \put(20,0){\circle{1}}
      \put(19.5,0.5){\line(-1,1){1}}

      \put(-0.2,-0.1){\(K_1\)}
      \put(1.7,1.8){\(K_2\)}
      \put(3.7,-0.1){\(K_3\)}
      \put(5.7,1.8){\(K_4\)}
      \put(7.7,-0.1){\(K_5\)}
      \put(9.7,1.8){\(K_6\)}
      \put(11.7,-0.1){\(K_7\)}
      \put(13.7,1.8){\(K_8\)}
      
      \put(17.4,2.9){\(K_{r-1}\)}
      \put(19.7,-0.1){\(K_r\)}
    \end{picture}
    }
      \caption{\(r\) odd}
            \label{fig:rodd}
    \end{figure}
    
    Hence, we greatest value is obtained for some \(K_i\) with \(i\) even.
    Furthermore, 
    \begin{equation}\label{eq:diff}
      K_{2i+2} - K_{2i} = c^{2i+1} \binom{r}{2i+1} +
      c^{2i+2}\binom{r}{2i+2} = c^{2i+1} \left( \binom{r}{2i+1} + c
        \binom{r}{2i+2} \right),
    \end{equation}
    which is \(>0\) iff
    \begin{equation}
      \label{eq:g0}
      c < - \frac{\binom{r}{2i+1}}{\binom{r}{2i+2}} = - \frac{2i+2}{r-2i-1}
    \end{equation}
    Since 
    \begin{equation}
      \label{eq:asc}
      - \frac{\binom{r}{3}}{\binom{r}{4}} > -
      \frac{\binom{r}{5}}{\binom{r}{5}} > \cdots > -
      \frac{\binom{r}{r-1}}{\binom{r}{r}} 
    \end{equation}
    the result for the odd case follows. The even case is proved
    similarly; here the inequalities for \(c < 0\) are as in
    Figure~\ref{fig:reven}. 
\begin{figure}[htbp]
  \tiny{
    \begin{picture}(22,3)
      \multiput(0,0)(4,0){4}{\circle{1}}
      \multiput(2,2)(4,0){4}{\circle{1}}
      \multiput(0.5,0.5)(4,0){4}{\line(1,1){1}}
      \multiput(3.5,0.5)(4,0){4}{\line(-1,1){1}}

      \put(18,0){\circle{1}}
      \put(20,2){\circle{1}}
      \put(18.5,0.5){\line(1,1){1}}

      \put(-0.2,-0.1){\(K_1\)}
      \put(1.7,1.8){\(K_2\)}
      \put(3.7,-0.1){\(K_3\)}
      \put(5.7,1.8){\(K_4\)}
      \put(7.7,-0.1){\(K_5\)}
      \put(9.7,1.8){\(K_6\)}
      \put(11.7,-0.1){\(K_7\)}
      \put(13.7,1.8){\(K_8\)}
      
      \put(17.4,1.0){\(K_{r-1}\)}
      \put(19.7,1.8){\(K_r\)}
    \end{picture}
    }
      \caption{\(r\) even}
      \label{fig:reven}
    \end{figure}
  \end{proof}
\end{section}

\begin{section}{Maximizing \(g\) on \(S\)}
  As we noted at the start of the previous section, if \(S\) consists
  of all square-free products of a finite set \(\set{p_1,\dots,p_r}\)
  of primes, then  \(\Delta(S)\) is an \((r-1)\)-simplex. Hence, if
  \(g\) is real-valued and \(g(s) \ge 1\)  (we call such a \(g\)
  \emph{multiplicative and 
    log-positive}), then the maximum of \(g\) on \(S\) is 
  \(g(p_1p_2 \cdots  p_r)\). More generally: 
  \begin{lemma}\label{lemma:maxfacet}
    Suppose that \(g\) is multiplicative and log-positive.
    Let \(W_1,\dots,W_\ell\) be the
    facets (i.e. a simplexes maximal w.r.t
    inclusion) of \(\Delta(S)\).
    Then the maximum value \(g(\sigma)\) for \(\sigma \in \Delta(S)\) is
    obtained on some facet \(W_i\). 
    
    Conversely, there exists a
    multiplicative and log-positive \(h\) so that
    \(h(W_1)\) is maximal.
  \end{lemma}
  \begin{proof}
    If \(\sigma \subset \tau\) then \(g(\sigma) \le g(\tau)\), so the
    maximum is attained on a facet.

    For the converse, define \(h\) on \(X(S)\) by 
    \begin{equation}
      \label{eq:hXS}
      h(p) = 
      \begin{cases}
        1 & \text{ if } p \notin W_1 \\
        2 & \text{ if } p \in W_1 
      \end{cases}
    \end{equation}
    We extend \(h\) to a multiplicative function on \(\Delta(S)\). It
    is then clear that \(h(W_1)=2^{\lvert W_1\rvert}\) whereas
    \(h(W_i)=2^{\lvert W_1 \cap W_i \rvert} < 2^{\lvert W_1\rvert}\)
    for \(i > 1\); the last inequality follows since \(W_1,W_i\) are
    facets and hence maximal w.r.t. inclusion.
    If we want a multiplicative \(h\) 
    which is strictly \(>1\) on non-empty simplexes, we can define
    \(h(p)=1+ \varepsilon\) for \(p \notin W_1\), where
    \(\varepsilon\) is some small positive number.
  \end{proof}

  We let \([n]=\set{1,2,\dots,n}\). Then \([n]\) is a unitary ideal,
  so we have
  \begin{corollary}
    If \(g\) is  multiplicative and log-positive function
    then the maximum \(g(s)\) with \(1 \le s \le n\) is obtained on a
    facet of \(\Delta([n])\).
  \end{corollary}
  
  As an example, if \(n=30\) then the 
  \(\Delta([30])\) looks like Figure~\ref{fig:d30}, so the 
  facets are 
  \[12, 14, 16, 17, 18, 19, 20, 21, 22, 23, 24, 25, 26, 27, 28, 29,
  30.\]
  Thus about \(57\%\) of the simplicies in \(\Delta([30])\) are
  facets.
  \begin{figure}[tbp]
    \centering
    \includegraphics[scale=0.8]{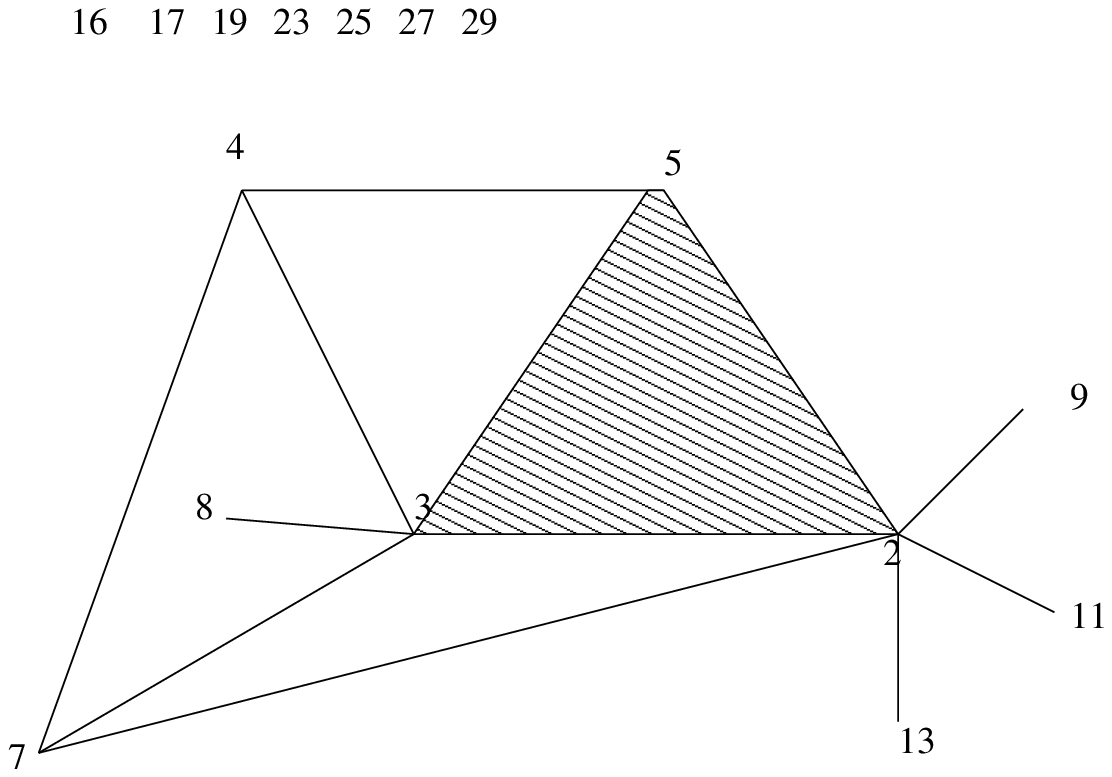}
    \caption{\(\Delta([30])\)}
    \label{fig:d30}
  \end{figure}
  In general, we have the following: 
  \begin{theorem}[Snellman]
    Let \(p_i\) denote the \(i\)'th prime number.
    For large \(n\), the number of facets in \(\Delta([n])\) is
    approximatively \(\gamma 
    n\), where 
    \begin{equation}
      \label{eq:7}
      \gamma = 1 - \frac{1}{2} + \sum_{i=1}^\infty \frac{ \frac{1}{p_i} -
    \frac{1}{p_{i+1}} }{\prod_{j=1}^i p_j} 
  \approx  0.607714359516618
    \end{equation}
  \end{theorem}
  \begin{proof}
    See \cite{Snellman:UninTrunc}.
  \end{proof}

  So, if we are to maximize (on \([n]\)) a large number of different
  \(g_i\)'s 
  which are  multiplicative and log-positive, it makes sense to
  precompute the facets of \(\Delta([n])\), and their
  factorizations. If \(q_1,\dots,q_r\) are the prime powers \(\le n\),
  and \(w_1,\dots,w_\ell\) are the facets of \(\Delta([n])\),
  let \(A=(a_{ij})\) be the \(\ell \times r\) integer matrix defined
  by
  \begin{equation}
    \label{eq:A}
    \forall 1 \le i \le \ell: \quad w_i= \prod_{j=1}^r q_j^{a_{ij}}
  \end{equation}
  Then, if  \(g\) is a log-positive multiplicative function,
  \begin{equation}
    \label{eq:Eval}
    \forall 1 \le i \le \ell: \quad g(w_i)= \prod_{j=1}^r g(q_j)^{a_{ij}}
  \end{equation}
  This means that in order to find the maximum for \(g\) we need to
  perform \(r\) evaluations to find the \(g(q_i)\)'s, then calculate
  \(\ell \approx \gamma n\) numbers, each of which is the product of
  at most
  \(r\) terms, and then find the maximum of those numbers. 
\end{section}

\begin{section}{Total orderings on \(S\) induced by \(g\)}
  As previously noted, if \(g\) is log-positive and \(\sigma \subset
  \tau\), then \(g(\sigma) \le g(\tau)\). Moreover, if \(g\) is
  \emph{strictly log-positive}, so that \(g(\sigma)>1\) for \(\sigma
  \neq \emptyset\), then 
  \begin{equation}
    \label{eq:strict}
    \sigma \subsetneq \tau \quad \implies \quad g(\sigma) < g(\tau).
  \end{equation}
  Assume that \(g\) has this property, that \(S\) is a unitary ideal,
  and that furthermore \(g\) is injective when restricted to \(S\).
  Then \(S\), and hence \(\Delta(S)\), is totally ordered by 
  \begin{equation}
    \label{eq:order}
    x > y \quad \iff \quad g(x) > g(y)
  \end{equation}
  It is clear, by \eqref{eq:strict}, that such a total order on
  \(\Delta(S)\) is a linear extension of the partial order given by
  inclusion of subsets. However, not all such linear extensions may
  occur.

  \begin{definition}
    Let \(r\) be a positive integer, and let \(V=\set{v_1,\dots,v_r}\)
    be a linearly ordered set with \(r\) 
    elements. Following MacLagan \cite{BoolTerm} we call a total order
    \(\succ\) on \(2^V\) a \emph{boolean termorder} if

      \begin{align}
      \label{eq:bt}
        \emptyset &\prec \sigma  & \text{ if } \emptyset \neq \sigma
        \subset V \\
        \sigma \cup \gamma & \prec \tau \cup \gamma & \text{ if } 
        \sigma \prec \tau \quad \text{ and } \quad \gamma \cap (\sigma
        \cup \tau) = \emptyset
      \end{align}
      We say that \(\prec\) is \emph{sorted} if 
      \begin{equation}
        \label{eq:sorted}
        v_1 \prec v_2 \prec \cdots \prec v_r
      \end{equation}
      Furthermore, \(\prec\) is \emph{coherent} if
      there exist \(r\) positive integers \(w_1,\dots,w_r\) such that
      \begin{equation}
        \label{eq:coher}
        \alpha \prec \beta \quad \iff \quad \sum_{v_i \in \alpha} w_i <
        \sum_{v_j \in \beta} w_j.
      \end{equation}
  \end{definition}
  \begin{lemma}\label{lemma:boolean}
    Suppose that \(S\) is a finite unitary ideal, and  let \(r\) be
    the number of prime powers in \(S\). Label 
    these prime powers \(v_1,\dots,v_r\). 
    Consider the set \(\mathcal{M}\) of  all multiplicative \(g\) 
    that are strictly
    log-positive,  injective when restricted to
    \(S\), and let \(\mathcal{M}^s\) denote the subset of those \(g\) that in
    addition fulfills
    \begin{equation}
      \label{eq:ord}
      g(v_1) < g(v_2) < \cdots < g(v_r).
    \end{equation}
    Let \(Y\) be the partial order on \(2^{\set{v_1,\dots,v_r}}\)
    which is generated by the  following relations: 
    \begin{equation}
      \label{eq:cover1}
        \set{v_{i_1},\dots,v_{i_\ell}}  <\!\!\!\cdot
        \set{v_{i_1},\dots,v_{i_\ell}, v_k} \text{ if } k \notin
        \set{i_1,\dots,i_\ell},
    \end{equation}
    \begin{multline}
        \set{v_{i_1},\dots,v_{i_j}, \dots, v_{i_\ell}}  <\!\!\!\cdot
        \set{v_{i_1},\dots,v_{i_j+1}, \dots, v_{i_\ell}} \\
        \text{ if } i_j +1
        \in [r] \setminus \set{v_{i_1},\dots,v_{i_j}, \dots, v_{i_\ell}}      
    \end{multline}
    Let \(T \subset S\), and let
     \(Y_T\) be the induced subposet on \(T \subset \Delta(S) \subset
    2^{\set{v_1,\dots,v_r}}\). Then 
    \begin{enumerate}[(i)]
    \item   Any total order on
      \(T\) induced by a \(g \in \mathcal{M}\) (by a \(g \in
      \mathcal{M}^s\)) is the 
      restriction of a  (sorted) coherent  boolean termorder.
    \item Conversely, the restriction to \(T\) of a (sorted) coherent boolean
    termorder 
    on \(2^{\set{v_1,\dots,v_r}}\) is induced by some \(g \in \mathcal{M}\) 
    (\(g  \in \mathcal{M}^s\). 
    \item   Any total order on
      \(T\) induced by a \(g \in \mathcal{M}^s\) is a linear
      extensions of \(Y_T\).
    \end{enumerate}
  \end{lemma}
  \begin{proof}
    We can W.L.O.G. assume that \(T=2^{\set{v_1,\dots,v_r}}\). If
    \(\prec\) is induced by \(g \in \mathcal{M}\) then 
    \begin{displaymath}
      \alpha \prec \beta 
      \quad \iff 
      \prod_{v_i \in \alpha}  g(v_i) < \prod_{v_j \in \beta} g(v_j) 
      \quad \iff 
      \sum_{v_i \in \alpha}  \log g(v_i) < \sum_{v_j \in \beta} \log g(v_j). 
    \end{displaymath}
    We can replace the \(\log g(v_i)\)'s by positive rational numbers
    that closely approximate them, and then, by multiplying out by a
    common denominator, by positive integers. Thus \(\prec\) is a
    coherent boolean termorder. It is clear that if \(g \in \mathcal{M}^s\) then
    \(\prec\) is sorted.
    
    If on the other hand \(\prec\) is a coherent boolean term order on
    \(2^{\set{v_1,\dots,v_r}}\), then there are positive integers
    \(w_1,\dots,w_r\) such that 
    \begin{displaymath}
      \alpha \prec \beta \quad \iff \sum_{v_i \in \alpha} w_i
      < \sum_{v_j \in \beta} w_j \quad \iff
      \prod_{v_i \in \alpha}  \exp(w_i) < \prod_{v_j \in \beta} \exp(w_j).
    \end{displaymath}
    If we define \(g(v_i)=w_j\) and extend this multiplicatively, then
    \(g\) induces \(\prec\).
    If \(\prec\) is sorted, clearly \(w_1 > w_2 > \cdots > w_r\), so
    \(g \in \mathcal{M}^s\).

    If \(g\) is strictly log-positive and fulfills \eqref{eq:ord} then
    clearly 
    \begin{equation}
      \label{eq:ordcons1}
      g(\set{v_{i_1},\dots,v_{i_\ell}}) < g(\set{v_{i_1},\dots,v_{i_\ell}, v_k})
    \end{equation}
    for \(k \notin \set{i_1,\dots,i_\ell}\), and likewise
    \begin{equation}
      \label{eq:ordcons2}
      g(\set{v_{i_1},\dots,v_{i_j}, \dots, v_{i_\ell}}) < 
      g(\set{v_{i_1},\dots,v_{i_j+1}, \dots, v_{i_\ell}})
    \end{equation}
    for 
    \(i_j +1  \in [r] \setminus \set{v_{i_1},\dots,v_{i_j}, \dots,
    v_{i_\ell}}.\) Thus any total order induced by a \(g \in \mathcal{M}\) is a
    linear extension of \(Y\).
  \end{proof}

  The symmetric group \(S_r\) acts transitively on 
  \(\mathcal{M}\), and  \(\setsuchas{\pi(\mathcal{M}^s)}{\pi \in
    S_r}\) is a partition of \(\mathcal{M}\) into \(\lvert S_r \rvert
  = r!\) blocks.
  Hence 
    \begin{corollary}\label{cor:nord}
      Let \(t(T)\) denote the number of total orders on \(T \subset
      \Delta(S)\) that are induced by multiplicative functions \(g \in
      \mathcal{M}^s\), and let \(\ell(Y_T)\) denote the number of linear
      extensions of \(Y_T\). Then 
      \begin{equation}
        \label{eq:ext}
        t(T) \le r! \ell(Y_T)
      \end{equation}
    \end{corollary}

    In the following example, we show
    that although any facet of \(\Delta(S)\) is maximal w.r.t. some
    total order induced by a log-positive multiplicative \(g\)
    (Lemma~\ref{lemma:maxfacet}), there are only certain orderings
    among those facets that are possible.
    
  \begin{example}
    Let us consider the poset \(Y\) with \(r=4\), and in particular
    the induced poset on the 2-subsets, which looks like
    Figure~\ref{fig:2pos}. 

    \newcommand{\goodgap}{\hspace{\subfigtopskip}\hspace{\subfigbottomskip}}
    \begin{figure}[hbtp]
      \caption{The poset \(Y\)}
    \setlength{\unitlength}{0.7cm}
    \centering
    \subfigure[The whole poset \(Y\)]{
    \includegraphics[scale=0.6]{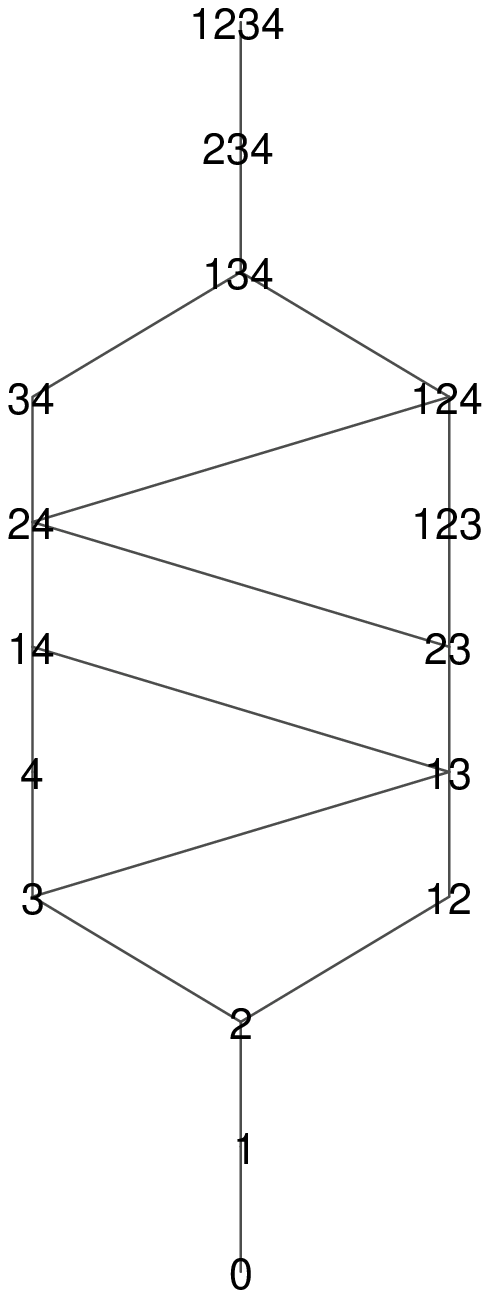}
    \label{fig:Y}
    }
    \hskip1cm
    \subfigure[\(Y\) restricted to 2-subsets]{
      \tiny{
    \begin{picture}(8,10)
      \put(3,1){\circle{1}}
      \put(3,3){\circle{1}}
      \put(1,5){\circle{1}}
      \put(5,5){\circle{1}}
      \put(3,7){\circle{1}}
      \put(3,9){\circle{1}}
      \put(3,1.5){\line(0,1){1}}
      \put(2.5,3.5){\line(-1,1){1}}
      \put(3.5,3.5){\line(1,1){1}}
      \put(1.5,5.5){\line(1,1){1}}
      \put(4.5,5.5){\line(-1,1){1}}
      \put(3,7.5){\line(0,1){1}}
      \put(2.7,0.9){12}
      \put(2.7,2.9){13}
      \put(0.7,4.9){23}
      \put(4.7,4.9){14}
      \put(2.7,6.9){24}
      \put(2.7,8.9){34}
    \end{picture}
  }
    \label{fig:2pos}
    }
    \end{figure}

    We see that this poset has exactly 2 linear
    extensions, corresponding to the two different ways of ordering
    the antichain \(\set{23,14}\). Thus, if 
    \begin{equation}
      \label{eq:ord4}
    g(v_1) < g(v_2) < g(v_3)  < g(v_4)
    \end{equation}
 then there are two possible orderings for
    \[\set{g(v_1v_2),g(v_1v_3), g(v_1v_4), g(v_2v_3), g(v_2v_4),
    g(v_3v_4)}.\]
  If we remove the restriction \eqref{eq:ord4}, then there are \(2
  \times 4! = 48\) different orderings, out of the \(\binom{4}{2}! =
  720\) \emph{a priori} possibilities. For instance, 
  \begin{equation}
    \label{eq:imposs}
    g(v_1v_2) >  g(v_2 v_4) > g(v_1v_3) > g(v_2v_3) > g(v_1v_4) >
    g(v_3v_4) 
  \end{equation}
  is impossible, since 
  \(g(v_2v_4) > g(v_1v_4) \implies g(v_2) > g(v_1)\) 
  but
  \(g(v_1v_3) > g(v_2v_3) \implies g(v_1) > g(v_2)\).
  Hence, there is no multiplicative arithmetic function \(g\) such
  that
  \begin{equation}
    \label{eq:imposs2}
    g(6) >  g(21) > g(10) > g(15) > g(14) > g(35).
  \end{equation}

  The whole poset \(Y\) on 4 letters looks like Figure~\ref{fig:Y}. It
  consists of \(16\) element an 
  has, as the reader may easily verify, 78 linear extensions.
  \end{example}
\end{section}

\bibliographystyle{hplain}
\bibliography{journals,articles,snellman}
\end{document}